\documentclass[a4paper, BCOR=0.0mm, DIV=calc]{scrartcl}
\usepackage[T1]{fontenc}
\usepackage[utf8]{inputenc}
\usepackage[ngerman]{babel}
\usepackage[osf,sc]{mathpazo}
\usepackage{textcomp}
\usepackage{microtype}
\usepackage{graphicx, color}
\usepackage{amsmath, amssymb, amsfonts}
\usepackage{booktabs}
\usepackage{nicefrac}
\usepackage{wasysym,pifont}
\usepackage{tikz}
\usepackage{graphicx}
\usepackage{float}
\usepackage{wrapfig}
\usetikzlibrary{calc,intersections}
\KOMAoptions{twoside=false, twocolumn=false, headinclude=false, footinclude=false, mpinclude=false, pagesize=auto}
\recalctypearea

\setkomafont{section}{\mdseries\scshape\Large}
\addtokomafont{subsection}{\mdseries\scshape}
\setkomafont{title}{\mdseries\scshape}

\begin{document}
\title{Bourlet's Theorem for the product of differential operators, an application of the operator method and a proof for $\sum_{n=1}^{\infty}\frac{1}{n^2}=\frac{\pi^2}{6}$, that Euler missed, derived from difference equations}
\author{Alexander Aycock}
\date{ }
\maketitle

\section*{Abstract}

We give another proof for

\[
\sum_{n=1}^{\infty}\frac{1}{n^2}=\frac{\pi^2}{6}
\]
that basically follows from the theory of difference equations.

\section*{Keywords}

Basel-Problem, Zeta-function values; Operator Calculus; Bourlet-Theory; infinite differential equations; differnce equations

\section*{Mathematics subject classification}

11M06, 39A06, 34A35; 

\section*{Introduction}

\textsc{Euler} got world-famous, when he proved in $1740$

\[
\sum_{n=1}^{\infty}\frac{1}{n^2}=\frac{\pi^2}{6}
\]
and solved the so-called Basel-Problem. In his career he gave at least six different proofs $[E41],[E61], [E63], [E189], [E464], [E592]$ for this series, all derived from completely different sources. In this paper we want to add another one, that could have also given by \textsc{Euler} himself.
This will also give us the oppurtunity to give an example of the operator calculus, which in its foundations traces back to \textsc{Heaviside} \cite{H} and is based on the formal inversion of a differential operator. We will begin with a proof of \textsc{Bourlet's} Theorem \cite {B}, which is important in the theory of differential equations of infinite order and that we will use, to solve the basic difference equation

\[
f(x+1)-f(x)=g(x)
\]
An application of the results will lead to another proof of the desired sum.

\section*{Derivation of Bourlet's Theorem for the product of differential equations}

\subsection*{Lemma 1 Leibniz-Rule}

Let $u(x)$ and $v(x)$ be functions $\in C^{\infty}(\mathbb{R}, \mathbb{R})$, then

\[
\frac{d^n}{dx^n}uv=\sum_{k=0}^{n} \binom{n}{k}u^{(n)}v^{(n-k)}
\]
This rule is well-known, that we can omit its proof here.

\subsection*{Lemma 2}

Let $F(x,z)=\sum_{n=0}^{\infty}F_n(x)z^n$, then formally

\[
\frac{F_z^{(k)}(x,z)}{k!}=\sum_{n=0}^{\infty}\binom{n+k}{k}F_{n+k}(x)z^n
\]
where $F_z^{(k)}(x,z)=\frac{\partial^k}{\partial z^k}F(x,z)$\\

$\mathbf{Proof}:$ We have

\begin{alignat*}{9}
&F_z^{k}(x,z) &&=\sum_{n=k}^{\infty} n(n-1) \cdots (n-k+1)F_n(x)z^{n-k} \\
& &&=\sum_{n=0}^{\infty} (n+k)(n+k-1) \cdots (n+1) F_{n+k}(x)z^n
\intertext{therefore}
&\frac{F_z^{(k)}(x,z)}{k!}&&=\sum_{n=0}^{\infty} \frac{n+k}{1} \cdot \frac{n+k-1}{2} \cdots \frac{n+1}{k}F_{n+k}(x)z^n\\
& &&=\sum_{n=0}^{\infty} \binom{n+k}{k}F_{n+k}(x)z^n
\end{alignat*}
which proves the resulz. Q.E.D.

\subsection*{Lemma 3}

Let $F(x,z)=\sum_{n=0}^{\infty} F_n(x)z^n$, and $u,v \in C^{\infty}(\mathbb{R}, \mathbb{R})$ and let $z=\frac{d}{dx}$, that $F(x,z)$ denotes a differential operator of infinite order, then formally

\[
F(x,z) uv= \sum_{n=0}^{\infty} \frac{v^{(n)}F_z^{(n)}}{n!}u
\]
$\mathbf{Proof}:$ Just use the $\textsc{Leibniz}-$rule for every term, then one finds

\begin{alignat*}{9}
&F_0(x)z^0~uv &&=F_0(x) &&\bigg[\binom{0}{0}u^{(0)}v^{(0)}\bigg]\\
&F_1(x)z^1~uv &&=F_1(x) &&\bigg[\binom{1}{0}u^{(1)}v^{(0)}&&+\binom{0}{1}u^{(0)}v^{(1)}\bigg]\\
&F_2(x)z^2~uv &&=F_2(x) &&\bigg[\binom{2}{0}u^{(2)}v^{(0)}&&+\binom{2}{1}u^{(1)}v^{(1)}&&+\binom{2}{1}u^{(1)}v^{(1)}\bigg]\\
&F_3(x)z^2~uv &&=F_3(x) &&\bigg[\binom{3}{0}u^{(3)}v^{(0)}&&+\binom{3}{1}u^{(2)}v^{(1)}&&+\binom{3}{1}u^{(1)}v^{(2)}&&+\binom{3}{3}u^{(3)}v^{(0)}\bigg]\\
& && && &&\text{etc.}
\end{alignat*}
Now just sum the columns

\begin{alignat*}{9}
&F(x,z)~uv &&= \sum_{k=0}^{\infty}v^{(k)}\sum_{n=0}^{\infty}F_{n+k}(x)\binom{n+k}{k}u^{(n)}\\
& &&=\sum_{n=0}^{\infty} \sum_{k=0}^{\infty} F_{n+k}(x)\binom{n+k}{k}v^{(k]}u^{(n)}\\
& &&=\sum_{n=0}^{\infty} \sum_{k=0}^{\infty} F_{n+k}(x)\binom{n+k}{k}v^{(k]}z^n ~u\\
& &&=\sum_{n=0}^{\infty} \frac{v^{(n)}F_z^{(n)}(x,z)}{n!}~u
\end{alignat*}
where Lemma $2$ was used in the last step. Q.E.D.

\subsection*{Bourlet's Theorem for the products of operators}

Let $X(x,z)=\sum_{n=0}^{\infty} X_n(x)z^n$ and $F(x,z)=\sum_{n=0}^{\infty}F_n(x)z^n$, then we have formally

\[
X(x,z)F(x,z)=\sum_{n=0}^{\infty} \frac{1}{n!}\frac{\partial^n X}{\partial z^n}\frac{\partial^n F}{\partial x^n}
\]
$\mathbf{Proof}:$ Consider the test function $u \in C^{\infty}(\mathbb{R}, \mathbb{R})$. By Lemma $3$ we have

\begin{alignat*}{9}
&X(x,z)F(x,z)~u&&= \sum_{j=0}^{\infty} \sum_{n=0}^{\infty} \frac{1}{n!}\frac{\partial ^n F_j}{\partial x^n}\frac{\partial^n X}{\partial z^n} ~u^{(j)}\\
& &&= \sum_{j=0}^{\infty} \sum_{n=0}^{\infty} \frac{1}{n!}\frac{\partial ^n F_j}{\partial x^n}\frac{\partial^n X}{\partial z^n} z^j ~u\\
& &&=\sum_{n=0}^{\infty}\frac{1}{n!}\frac{\partial^n X}{\partial z^n}\sum_{j=0}^{\infty} \frac{\partial^n F_j}{\partial x^n}z^j ~u\\
& &&=\sum_{n=0}^{\infty}\frac{\partial^n X}{\partial z^n}\frac{\partial^n F}{\partial x^n}~u
\end{alignat*}
which proves the theorem. Q.E.D.\\

\textsc{Bourlet's} Theorem is crucial for the study of differential equations of infinie order and seems to be not very well-known despite its importance. Hence it was appropriate, to mention and derive it here. We will just need the following corollary.

\subsection*{Corollary}

If $F$ is independend of $x$, then the iverse operator is easily found, because we just have to solve the equation

\[
X(x,z)F(x,z)=\sum_{n=0}^{\infty} \frac{1}{n!}\frac{\partial^n X}{\partial z^n}\frac{\partial^n F}{\partial x^n}=1
\]
and because $F$ is independend of $x$ wie zu have

\begin{alignat*}{9}
&X(z)F(z)&&=1 \quad \text{and}\\
&X(z) &&=\frac{1}{F(z)}
\end{alignat*}
Now we can use this result to solve the basic difference equation.

\section*{Solution of the difference equation $f(x+1)-f(x)=g(x)$}

At first we note, that we have by \textsc{Taylor's} Theorem

\[
f(x+a)=\sum_{n=0}^{\infty} \frac{f^{(n)}(x)}{n!}a^n
\]
and hence for $a=1$
\[
f(x+1)=\sum_{n=0}^{\infty} \frac{f^{(n)}(x)}{n!}
\]
So if we substitute this into the propounded difference equation we obtain

\begin{alignat*}{9}
& \sum_{n=1}^{\infty} \frac{f^{(n)}(x)}{n!}&&=g(x)\\
& \sum_{n=1}^{\infty} \frac{z^n}{n!}f(x)&&=g(x)\\
&\quad (e^z-1) &&=g(x)
\end{alignat*}
We will need the following partial fraction decomposition of $\frac{1}{e^z-1}$

\[
\frac{1}{e^z-1}=-\frac{1}{2}+\frac{1}{z}+ \sum_{n=1}^{\infty} \frac{1}{z-2n \pi i}+\frac{1}{z+2n \pi i}
\]
which is easily proved by calculating the residues and an application of \textsc{Liouville's} Theorem. \\
Now we just proceed in the usual way, to  solve a differential equation, and solve the homogenous equation at first, namely

\[
f(x+1)=f(x)
\]
which is easily seen to be solved in general by

\[
f_{hom}(x)=\sum_{n=-\infty}^{\infty} c_n e^{2 \pi ni x}
\]
a Fourier series, of course. To find one particular solution, we just apply the operator

\[
\frac{1}{e^z-1}=-\frac{1}{2}+\frac{1}{z}+ \sum_{n=1}^{\infty} \frac{1}{z-2n \pi i}+\frac{1}{z+2n \pi i}
\]
to $g(x)$, note, that $\frac{1}{z}$ is to be interpreted as $\int^x_a$, meaning the integral. Therefore we find

\begin{alignat*}{9}
&f_{part}(x)&&=\bigg(-\frac{1}{2}&&+~~~~~~\frac{1}{z}&&+ \sum_{n=1}^{\infty} ~~~~~~\frac{1}{z-2n \pi i}&&+~~~~~~\frac{1}{z+2n \pi i} \bigg)g(x)\\
& &&=-\frac{1}{2}g(x)&&+\int^x g(t)dt&&+ \sum_{n=1}^{\infty} \sum_{k=0}^{\infty} (2n \pi i)^k \int^{k+1}g(x)dx&&+(-2n \pi i)^k \int ^{k+1}g(x)dx\\
& &&=-\frac{1}{2}g(x)&&+\int^x g(t)dt&&+ \sum_{n=1}^{\infty} \sum_{k=0}^{\infty} (2n \pi i)^k \int^{x}\frac{(x-t)^kg(t)dt}{k!}&&+(-2n \pi i)^k \int ^{x}\frac{(x-t)^kg(t)dt}{k!}\\
& &&=-\frac{1}{2}g(x)&&+\int^x g(t)dt&&+ \sum_{n=1}^{\infty}   \int^{x}e^{2n \pi i(x-t)}g(t)dt&&+\int ^{x}e^{-2n \pi i(x-t)}g(t)dt\\
\end{alignat*}
where we used the well-known formula to reduce the iterated integral to a single one, namely

\[
\int ^{k+1}f(x)dx=\int^x \frac{(x-t)^kf(t)}{k!}dt
\]
Note, that for $g(x)=0$ we would obtain the Fourier series again. Our general  solution, if $\Pi(x+1)=\Pi(x)$ is a periodic function

\[
f(x)=\Pi(x)-\frac{1}{2}g(x) \int ^x g(t)dt +2 \sum_{n=1}^{\infty} \int^x \cos(2n \pi(x-t))g(t)dt
\]
where we omitted the lower boundaries in the integrals for the sake of brevity.

\section*{Proof of the series $\sum_{n=1}^{\infty} \frac{1}{n^2}=\frac{\pi^2}{6}$}

Now we can just apply our results to prove the propounded series

Consider the sum

\[
\sum_{k=1}^{x} (k-1)^2=f(x)
\]
then one by elementary methods one finds fir integral $x$

\[
f(x)=\frac{x^3}{3}-\frac{x^2}{2}-\frac{x}{6}
\]
Now we can rewrite this as difference equation

\[
f(x+1)-f(x)=x^2
\]
and we already now, that the value found above is a particular solution, comparing it to our general solution we obtain

\begin{alignat*}{9}
&\frac{x^3}{3}&&-\frac{x^2}{2}&&+\frac{x}{6}&&=\Pi(x)&&-\frac{1}{2}x^2 &&+\int ^x x^2dt &&+2 \sum_{n=1}^{\infty} \int^x \cos(2n \pi(x-t))t^2dt\\
& && && &&=\Pi(x)&&-\frac{1}{2}x^2 &&+\frac{x^3}{3} &&+\frac{x}{\pi^2} \sum_{n=1}^{\infty} \frac{1}{n^2}
\end{alignat*}
where all constants of integration can be absorbed in the periodic function.

Now for $x=0$ we conclude

\[
0=\Pi(0)=\Pi(n) \quad \text{for all} \quad  n \in \mathbb{N}
\]
So let $n$ be a natural number, then

\[
\frac{n^3}{3}-\frac{n^2}{2}+\frac{n}{6}=\Pi(n)+\frac{n^3}{3}-\frac{n}{2}-\frac{n}{\pi^2}\sum_{k=1}^{\infty}\frac{1}{k^2}
\]
and by comparing he coefficients

\[
\frac{1}{6}=\frac{1}{\pi^2}\sum_{k=1}^{\infty}\frac{1}{k^2}
\]
and hence

\[
\sum_{k=1}^{\infty}\frac{1}{k^2}=\frac{\pi^2}{6}
\]
which is the desired series.

\section*{Concluding Remarks}

In similar manner we can conclude

\begin{alignat*}{9}
&\sum_{n=1}^{\infty} \frac{\pi^4}{90}&& \quad \text{from} \quad &&g(t)=t^4\\
&\sum_{n=1}^{\infty} \frac{\pi^6}{945}&& \quad \text{from} \quad &&g(t)=t^6\\
&\sum_{n=1}^{\infty} \frac{\pi^8}{9450}&& \quad \text{from} \quad &&g(t)=t^8\\
& && \quad \text{etc.}
\end{alignat*}
This proof could indeed have been given by \textsc{Euler}, who was again the first to consider differential equations of infinie orde in $[E189]$, but he considered other things in this paper and made a slight mistake in the solution of the corresponding differential equation.
We could proceed further and prove other sums of this kind, we would just have to find the corresponding difference equation and solve it by the operator method. But we will undertake this task on another occasion.


\begin{thebibliography}{99}
\bibitem[B]{B} C. Bourlet {\em Sur les operations en general et les equations differentilles lineaires d'ordre infini}, 	Annales de l'Ecole Normale Superieure, 3rd ser., vol. 14 (1897), pp.130-190
\bibitem[E41]{E41} L. Euler {\em De summis serierum reciprocarum}, 	 Opera Omnia: Series 1, Volume 14, pp. 73 - 86 
\bibitem[E61]{E61} L. Euler {\em De summis serierum reciprocarum ex potestatibus numerorum naturalium ortarum dissertatio altera, in qua eaedem summationes ex fonte maxime diverso derivantur }, 	 Opera Omnia: Series 1, Volume 14, pp. 138 - 155 
\bibitem[E63]{E63} L. Euler {\em Demonstration de la somme de cette suite $1+\frac{1}{4}+\frac{1}{9}+\frac{1}{16}+$ etc. }, 	 Opera Omnia: Series 1, Volume 14, pp. 177 - 186 
\bibitem[E189]{E189} L. Euler {\em De serierum determinatione seu nova methodus inveniendi terminos generales serierum }, 	  Opera Omnia: Series 1, Volume 14, pp. 463 - 515 
\bibitem[E464]{E464} L. Euler {\em Nova methodus quantitates integrales determinandi }, 	  Opera Omnia: Series 1, Volume 15, pp. 621 - 660
\bibitem[E592]{E594} L. Euler {\em De resolutione fractionum transcendentium in infinitas fractiones simplices },  Opera Omnia: Series 1, Volume 17, pp. 421 - 457 
\bibitem[H]{H} O. Heaviside {\em Electromagnetic Theory, Volume 1-3}, Forgotten Books, 2010
\end{thebibliography}
\end{document}